\documentclass[12pt]{amsart}
\usepackage{amscd,amssymb}
\usepackage[graph,frame,poly,arc]{xy}  
\usepackage[plainpages,backref,urlcolor=blue]{hyperref}

\theoremstyle{plain}
\newtheorem{thm}[subsection]{Theorem}

\newtheorem{prop}[subsection]{Proposition}
\newtheorem{cor}[subsection]{Corollary}

\theoremstyle{definition}
\newtheorem{rk}[subsection]{Remark}

\newtheorem{ex}[subsection]{Example}
\newtheorem{conj}[subsection]{Conjecture}

\numberwithin{equation}{section}
\newcommand{\I}{{\mathcal I}}

\newcommand{\A}{{\mathcal A}}
\newcommand{\B}{{\mathcal B}}

\newcommand{\al}{{\alpha}}

\newcommand{\C}{\mathbb{C}}
\newcommand{\PP}{\mathbb{P}}

\DeclareMathOperator{\im}{im}

\DeclareMathOperator{\indeg}{indeg}
\begin{document}
%\date{June 4, 2009}

\title [Plane curves with three syzygies]
{Plane curves with three syzygies, minimal Tjurina curves, and nearly cuspidal curves}

\author[Alexandru Dimca]{Alexandru Dimca$^{1}$}
\address{Universit\'e C\^ ote d'Azur, CNRS, LJAD and INRIA, France and Simion Stoilow Institute of Mathematics,
P.O. Box 1-764, RO-014700 Bucharest, Romania}
\email{dimca@unice.fr}

\author[Gabriel Sticlaru]{Gabriel Sticlaru}
\address{Faculty of Mathematics and Informatics,
Ovidius University
Bd. Mamaia 124, 900527 Constanta,
Romania}
\email{gabrielsticlaru@yahoo.com }

\thanks{$^1$ This work has been partially supported by the French government, through the $\rm UCA^{\rm JEDI}$ Investments in the Future project managed by the National Research Agency (ANR) with the reference number ANR-15-IDEX-01 and by the Romanian Ministry of Research and Innovation, CNCS - UEFISCDI, grant PN-III-P4-ID-PCE-2016-0030, within PNCDI III}



\begin{abstract} 
We start the study of reduced complex projective plane curves, whose Jacobian syzygy module has 3 generators. Among these curves one finds the nearly free curves introduced by the authors, and the plus-one generated line arrangements introduced by Takuro Abe. All the Thom-Sebastiani type plane curves, and more generally, any curve whose global Tjurina number is equal to a lower bound given by A. du Plessis and C.T.C. Wall, are 3-syzygy curves.
Rational plane curves which are nearly cuspidal, i.e. which have only cusps except one singularity with two branches, are also related to this class of curves.
\end{abstract}
 
\maketitle

%\tableofcontents

\section{Introduction} 

Let $S=\C[x,y,z]$ be the polynomial ring in three variables $x,y,z$ with complex coefficients, and let $C:f=0$ be a reduced curve of degree $d\geq 3$ in the complex projective plane $\PP^2$. 
We denote by $J_f$ the Jacobian ideal of $f$, i.e. the homogeneous ideal in $S$ spanned by the partial derivatives $f_x,f_y,f_z$ of $f$, and  by $M(f)=S/J_f$ the corresponding graded quotient ring, called the Jacobian (or Milnor) algebra of $f$.
Consider the graded $S$-module of Jacobian syzygies of $f$, namely
$$AR(f)=\{(a,b,c) \in S^3 \ : \ af_x+bf_y+cf_z=0\}.$$
According to Hilbert Syzygy Theorem, the graded Jacobian algebra $M(f)$ has a minimal resolution of the form
\begin{equation}
\label{res1}
0 \to F_3 \to F_2 \to F_1 \to F_0,
\end{equation}
where clearly $F_0=S$, $F_1=S^3(1-d)$ and the morphism $F_1 \to F_0$ is given by
$$(a,b,c) \mapsto af_x+bf_y+cf_z.$$
With this notation, the graded $S$-module of Jacobian syzygies $AR(f)$ has the following minimal resolution
$$0 \to F_3(d-1) \to F_2(d-1).$$
We say that $C:f=0$ is an {\it $m$-syzygy curve} if the module  $F_2$ has rank $m$. Then the module $AR(f)$ is generated by $m$ homogeneous syzygies, say $r_1,r_2,...,r_m$, of degrees $d_j=\deg r_j$ ordered such that $$1 \leq d_1\leq d_2 \leq ...\leq d_m.$$ 
We call these degrees the {\it exponents} of the curve $C$ and $r_1,...,r_m$ a {\it minimal set of generators } for the module  $AR(f)$. Note that $d_1=mdr(f)$ is the minimal degree of a non trivial Jacobian relation in $AR(f)$. It is known that, for a reduced degree $d$ curve $C$, one has $d_m \leq 2d-4,$
 see \cite[Corollary 11]{CD}, and for the case $C$ a line arrangement, one has the much stronger inequality
 $d_m \leq d-2,$
see \cite[Corollary 3.5]{Sch0}. On the other hand,
for a reduced degree $d$ curve $C$, one has 
$m \leq d+1,$
see \cite[Proposition 2.1]{DStMax}, while for the case $C$ a line arrangement, we have the slightly stronger inequality
 $m \leq d-1,$
see \cite[Corollary 1.3]{DStSat}, as well as its generalization in Corollary \ref{corB1} below. The curve $C$ is {\it free} when $m=2$, since then  $AR(f)$ is a free module of rank 2, see for such curves \cite{B+,Dmax,DStFD,Sim1,Sim2,ST}. In this note we consider the next simplest possibility for this resolution, namely we start the study of the 3-syzygy curves. {\it However, some of our main results, see Theorems
\ref{thmPO1}, \ref{thmd3}, \ref{thmNEW}, \ref{thmPO2}, concern any reduced plane curve.}

\begin{ex}
\label{ex1} Here are two classes of 3-syzygy curves which are already intensely studied.

\noindent (i) The {\it nearly free curves}, introduced in \cite{DStRIMS} and studied in \cite{AD, B+, Dmax, Drcc, MaVa} correspond to the 3-syzygy curves satisfying $d_3=d_2$ and $d_1+d_2=d$. 

\noindent (ii) The {\it plus-one generated line arrangements} of level $d_3$, introduced by Takuro Abe in \cite{Abe18},
corresponds to 3-syzygy line arrangements satisfying $d_1+d_2=d$. A 3-syzygy curve will be called a {\it plus-one generated curve}  if it satisfies $d_1+d_2=d$. In fact, in Theorem \ref{thmPO1} we show that the condition $d_1+d_2=d$ characterizes the plus-one generated curves among all the reduced plane curves. 
Both nearly free curves and the plus-one generated curves satisfy
$$mdr(f)=d_1 \leq \frac{d}{2}.$$
The line arrangement
$$ C: f=xyz(x+y-2z)(x-3y+z)(-5x+y+z)(x+y+z)=0$$
is a 3-sygyzy curve with exponents $d_1=d_2=d_3=4$ and  $\tau(C)=24$, but it is not a plus-one generated arrangement since $d_1+d_2>d=7$.

\end{ex}

Note that if $C: f=0$ be a {\it Thom-Sebastiani plane curve}, i.e. a curve such that $f(x,y,z)=g(x,y)+z^d$, where
$g$ is a homogeneous polynomial of degree $d$ in $R=\C[x,y]$, then $C$ is always a 3-syzygy curve. Moreover, these curves realize the equality in the {\it du Plessis-Wall inequality}
for the global Tjurina number $\tau(C)$ in \ref{dPW}, see Example \ref{ex1.1} below. Conversely, such a minimal Tjurina curve, realizing the equality in the du Plessis-Wall inequality, is a 3-syzygy curve satisfying $d_2=d_3=d-1$,  Theorem \ref{thmNEW}. This result is surprising, since the free curves correspond to $m=2$ and realize the maximal possible values for the global Tjurina number $\tau(C)$, see
\cite{Dmax,dPW}.

\medskip

Our main interest in the free and nearly free curves comes from the following.
\begin{conj}
\label{c1}
A reduced plane curve $C:f=0$ which is rational cuspidal is either free, or nearly free.
\end{conj}
This conjecture is known to hold when the degree of $C$ is even, as well as in many other cases, in particular for all odd degrees $d \leq 33$, see \cite{ Drcc, DStRIMS,DStMos}.
We call an irreducible  plane  curve, having only cusps (i.e. unibranch singularities), except for one singular point which has two branches, a {\it nearly cuspidal curve}. The simplest example of such a curve is a nodal cubic, for instance
$C:f=xyz+x^3+y^3=0$, for which the module $AR(f)$ is minimally generated by four syzygies, all of degree 2.  Other examples of rational, nearly cuspidal curves of odd degree which are 4-syzygy curves are given below in Example \ref{exNU=2}.

The following is one of our main results, obtained by combining Theorem \ref{thmPO2} and Theorem \ref{thm2} below.
\begin{thm}
\label{thmMAIN}
Let $C$ be a plane, rational, nearly cuspidal curve of even degree $d \geq 2$.
Then $C$ is a free curve, or a nearly free curve or a plus-one generated curve.
\end{thm}
Examples of plus-one generated curves of this type can be found in Example 
\ref{exlowdegree} (i), (ii) and (iii), and in Example \ref{ex2}. On the other hand, examples of rational nearly cuspidal curves which are free can be found in \cite{N}.

\medskip

In the second section we start with some basic facts on 3-syzygy curves, namely the formulas for the main invariants of $f$ expressed in terms of the exponents $(d_1,d_2,d_3)$ given in Proposition \ref{propA}. This leads to a characterization of plus-one generated curves among all plane curves in Theorem \ref{thmPO1}.
In the third section we prove that  $C$ is a 3-sygyzy curve if and only if the associated Bourbaki ideal $B(C,r_1)$, which is always a locally complete intersection, is in fact a global complete intersection, see Proposition \ref{propD}.
 This leads to a characterizations of plus-one generated curves among the 3-syzygy curves, see Proposition \ref{propB},
as well as a characterization of 3-syzygy curves, see Corollary \ref{corBour}, both of them in terms of the total Tjurina number $\tau(C)$, which is the sum of the Tjurina numbers of the singular points of $C$.
In Theorem \ref{thmNEW}, we also give a refined version of  the  du Plessis-Wall lower bound
for the global Tjurina number in \ref{dPW}.

Let $I_f$ denote the saturation of the ideal $J_f$ with respect to the maximal ideal ${\bf m}=(x,y,z)$ in $S$ and consider the following  local cohomology group, usually called the Jacobian module of $f$, 
 $$N(f)=I_f/J_f=H^0_{\bf m}(M(f)).$$
We set $n(f)_k=\dim N(f)_k$ for any integer $k$ and also $\nu(C)=\max _j \{n(f)_j\}$, as in \cite{AD,Drcc}.
Note that $C$ is free if and only if $\nu(C)=0$, and $C$ is nearly free if and only if $\nu(C)=1$.
In this note we show that $\nu(C)=2$ implies that $C$ is a plus-one generated curve when $d$ is even. When $d=2d'-1$ is odd, then  $C$ is either a plus-one generated curve, or $C$ is a 4-syzygy curve with exponents $d_1=d_2=d_3=d_4=d'$, see Theorem \ref{thmPO2}. In Example \ref{exNU=2} we give several examples of such 4-syzygy curves, and we show that for a line arrangement $\A$ with $\nu(\A)=2$, the numbers of $k$-multiple points for all $k$ do not determine whether $\A$ is a plus-one generated arrangement or not.

Our interest in curves with $\nu(C)=2$ is motivated by the following.
\begin{conj}
\label{c2}
A plane, rational, nearly cuspidal  curve $C:f=0$ satisfies 
$$\nu(C) \leq 2.$$
\end{conj}

In Theorem \ref{thmHS}  and Corollary \ref{corHS} we describe the resolution of the Jacobian module $N(f)$ for a 3-syzygy curve, and respectively the Hilbert function of $N(f)$ for a plus-one generated curve, in particular we show that $\nu(C)=d_3-d_2+1$, which may be arbitrarily large. Theorem \ref{thmPO2} says that this type of Hilbert function occurs exactly for plus-one generated curves. Hassanzadeh and Simis result in \cite[Proposition 1.3]{HS} plays a key role in these proofs.

In the fourth section we give various examples of plus-one generated curves and 3-syzygy curves: irreducible curves of low degree in Example \ref{exlowdegree}, curves obtained from a smooth curve plus a secant line in Example \ref{exConj}, and an infinite series of plus-one generated irreducible curves in Example \ref{ex2}. Example \ref{exNUlarge} shows that the invariant $\nu(C)$ can be arbitrarily large for a plus-one generated curve.

In the final section we give some upper bounds for the degrees of a minimal set of generators for the module $AR(f)$, in terms of genera of the irreducible components of $C:f=0$. We also prove  Conjecture \ref{c2} in Theorem \ref{thm2}  when the degree $d$ of $C$ is even, or $d=p^k$ with $p$ prime, or  when $d$ is odd and some extra topological condition is satisfied by the singularities of $C$.
Since this topological condition is hard to check in practice, we describe other situations where
Conjecture \ref{c2} holds in Theorem \ref{thm3}, Theorem \ref{thm4} and Corollary \ref{corConj2}, covering in particular all the odd degrees $d \leq 33$.

\medskip

We would like to thank Philippe Ellia for allowing us to use his idea for the proof of Theorem \ref{thmd3}, thus answering an open question in a preliminary version of this paper.

\section{Results involving the exponents}

With the above notation, for $C:f=0$ a 3-syzygy curve, the resolution \eqref{res1} becomes
\begin{equation}
\label{res2}
0 \to S(-e) \to \oplus_{j=1} ^3S(1-d-d_j)\to S^3(1-d)  \to S.
\end{equation}
Using \cite[Lemma 1.3]{HS}, we know that $e\geq d+d_3$. Moreover, for $k\geq e-2$, one has the obvious formula
\begin{equation}
\label{eq1}
m(f)_k=\dim M(f)_k={k+2 \choose 2}-3{k-d+3 \choose 2}+\sum_{j=1} ^3{k-d+3+d_j \choose 2}-{k-e+2 \choose 2}.
\end{equation}
Since one has $m(f)_k=\tau(C)$, the total Tjurina number of the curve $C$, for all large $k$, it follows that the polynomial in $k$ given by the alternated sum in \eqref{eq1} must be a constant. This implies our first result.

\begin{prop}
\label{propA}
Suppose $C:f=0$ is a 3-syzygy curve of degree $d\geq 3$ with exponents $(d_1,d_2,d_3)$. Then 
the following hold.
\begin{enumerate}

\item $e=d_1+d_2+d_3$. In particular $d_1+d_2 \geq d$.

\item $ct(f)=d-2+d_1$.

\item $st(f)=reg(f)+1=d_1+d_2+d_3-2$.

\item $\tau(C)=(d-1)(d_1+d_2+d_3)-(d_1d_2+d_1d_3+d_2d_3)$.

\end{enumerate}

\end{prop}
Here $ct(f)$ is the {\it coincidence threshold} of $f$, $st(f)$ is the {\it stability threshold} of $f$, $reg(f)$ is the {\it Castelnuovo-Mumford regularity} of the Milnor algebra $M(f)$, 
and $\tau(C)$ is the {\it total Tjurina number} of $C$, see for instance \cite{Dmax, DIM, DStFD,DStRIMS} for the corresponding definitions.
Note that $d_1+d_2 \leq d-1$ implies that $C$ is a free curve, see \cite{ST}.

\begin{cor}
\label{corA}
Suppose $C:f=0$ is a 3-syzygy curve of degree $d\geq 3$ with a minimal set of generators $r_1,r_2,r_3$ for $AR(f)$ and exponents $(d_1,d_2,d_3)$. Then the syzygies among $r_1,r_2,r_3$ are generated by a unique relation
$$R: h_1r_1+h_2r_2+h_3r_3=0,$$
where $h_i$ are homogeneous polynomials in $S$, with $\deg h_i=d_j+d_k-d+1\geq 1$ for any permutation $(i,j,k)$ of $(1,2,3)$. Moreover,
one has $g.c.d. (h_i,h_j)=1$ for any $i\ne j$. In particular $h_i \ne 0$ for all $i=1,2,3$.
\end{cor}
The last claim is related to the notion of strictly plus-one generated line arrangements, see \cite[Definition 1.1 and Proposition 5.1]{Abe18}. See also Remark \ref{rkcoeff} below for more information on $h_2$ and $h_3$.
\proof
The resolution \eqref{res2} and the formula for $e$ in Proposition \ref{propA} give the first claim.
To prove the remaining ones, suppose $i=1$, $j=2$ and that the irreducible homogeneous polynomial $c$ is a common factor of $h_1$ and $h_2$, with $\deg c >0$. Then there are two possibilities.
The first one is that $c$ divides $h_3$ as well. Then we can simplify in $R$ by $c$ and get a new relation
$$R': h'_1r_1+h'_2r_2+h'_3r_3=0,$$
of strictly lower degree than $R$, a contradiction.
If $c$ does not divide $h_3$, then it must divide all the components $(a_3,b_3,c_3)$ of $r_3$.
Dividing by $c$ these components, we get a new relation $r_3'=(a'_3,b'_3,c'_3) \in AR(f)$, of degree $<d_3$. This shows that $r_3$ is already in the submodule spanned by $r_1$ and $r_2$,
a contradiction. The other cases can be treated similarly.
\endproof

We give next a characterization of plus-one generated curves among all the reduced plane curves. 
Consider the general form of the minimal resolution for the Milnor algebra $M(f)$ of a curve $C:f=0$ that is assumed to be not free, namely
\begin{equation}
\label{res2A}
0 \to \oplus_{i=1} ^{m-2}S(-e_i) \to \oplus_{j=1} ^mS(1-d-d_j)\to S^3(1-d)  \to S,
\end{equation}
with $e_1\leq ...\leq e_{m-2}$ and $d_1\leq ...\leq d_m$.
It follows from \cite[Lemma 1.1]{HS} that one has
\begin{equation}
\label{res2B}
e_j=d+d_{j+2}-1+\epsilon_j,
\end{equation}
for $j=1,...,m-2$ and some integers $\epsilon_j\geq 1$. Using the same approach as for the proof of Proposition \ref{propA}, namely that one must have $m(f)_k=\tau(C)$ for all large $k$, or even better using \cite[Formula (13)]{HS}, it follows that one has
\begin{equation}
\label{res2C}
d_1+d_2=d-1+\sum_{i=1} ^{m-2}\epsilon_j.
\end{equation}
This formula implies the following.

\begin{thm}
\label{thmPO1}
Let $C:f=0$ be a reduced plane curve of degree $d$ and let $d_1$ and $d_2$ be the minimal degrees of a minimal system of generators for the module of Jacobian syzygies $AR(f)$ as above.
Then the following hold.
\begin{enumerate}

\item The curve $C$ is free if and only if $d_1+d_2=d-1$.

\item The curve $C$ is plus-one generated if and only if $d_1+d_2=d$.

\end{enumerate}

\end{thm}

\proof
The claim (1) is well known, see for instance \cite{ST}, and it is included only to point out the perfect analogy of the class of free curves and the class of plus-one generated curves.
The second claim follows from the formula \eqref{res2C} and the inequalities $\epsilon_j\geq 1$.

\endproof

We end this section with the following result.

\begin{thm}
\label{thmd3}

 Let $C:f=0$ be an $m$-syzygy curve of degree $d$ with exponents $1\leq d_1 \leq d_2\leq \cdots  \leq d_m$ and $m \geq 3$. Then $d_1+d_2 \geq d$ and 
 $$d_1 \leq d_2 \leq d_3 \leq d-1.$$

\end{thm}

\proof
Note that the condition $m \geq 3$ is equivalent to the fact that $C$ is not a free curve.
Then the claim $d_1+d_2 \geq d$ follows from \cite[Lemma 1.4]{ST}.

\bigskip

The proof of the claim $d_3 \leq d-1$, which is the main claim above,  is based on an {\it idea of Philippe Ellia}. First note that the claim does not depend on the choice of linear coordinates on $\PP^2$. So we choose these coordinates such that the line $L:x=0$ is transversal to $C$, i.e. such that the
binary form $f(0,y,z)$ has $d$ distinct linear factors. This implies that the ideal $(f_y,f_z)$ defines a 0-dimensional subscheme in $\PP^2$.

Let $r_1 \in AR(f)$ be a non-zero syzygy of minimal degree $d_1$, and consider the graded quotient module $\overline {AR(f)}=AR(f)/S r_1$. The exact sequence \eqref{B3} below implies that the $S$-module $\overline {AR(f)}$ is torsion free.

Assume that $d_3 \geq d$ and let $r_2 \in AR(f)$ be a syzygy of degree $d_2$, part of a minimal system of generators for $AR(f)$ starting with $r_1$.
Consider now the Koszul syzygies in $AR(f)_{d-1}$ given by
$$k^x=(0,f_z,-f_y) \ k^y=(-f_z,0,f_x) \text{ and } k^z=(f_y,-f_x,0).$$
 Then one can write
$$k^x=a^xr_1+b^xr_2, k^y=a^yr_1+b^yr_2 \text{ and } k^z=a^zr_1+b^zr_2, $$
for some polynomials $a^x,a^y,a^z \in S_{d-1-d_1}$ and $b^x,b^y,b^z \in S_{d-1-d_2}$.
The obvious relation
$$f_xk^x+f_yk^y+f_zk^z=0,$$
yields the relation
\begin{equation} \label{tors}
Ar_1+ Br_2=0,
\end{equation}
where $A=a^xf_x+a^yf_y+a^zf_z$ and $B=b^xf_x+b^yf_y+b^zf_z$.
There are two cases to discuss.

\medskip

\noindent {\bf Case 1. $B \ne 0$.} 
In this case we reach a contradiction, since $B \ne 0$ and equation \eqref{tors}
implies that the class of $r_2$ in $\overline {AR(f)}$ is a non-zero torsion element.

\medskip

\noindent {\bf Case 2. $B=0$.} 
Then the inequality $d_1+d_2 \geq d$  implies that
$$d-1-d_2\leq d_1-1,$$
and hence $b^x=b^y=b^z=0$, since $d_1$ is the minimal degree of a non trivial Jacobian syzygy.
This vanishing in turn implies that
$$f_y=a^za_1 \text{ and } f_z=-a^ya_1,$$
where $r_1=(a_1,b_1,c_1)$. This is a contradiction with the fact that the ideal $(f_y,f_z)$ defines a 0-dimensional subscheme in $\PP^2$.
\endproof

\begin{rk}
\label{rkS1}
 It is known that, for any singular curve $C:f=0$, one has $st(f) \leq T=3d-6$, see \cite{DBull}. It follows that, for a 3-syzygy singular curve $C:f=0$, one has
$$d_1+d_2+d_3 \leq 3d-4.$$
Note that, for a uninodal plane curve of degree $d$, the module $AR(f)$ has 4 minimal generators, of degrees $d-1, d-1, d-1, 2d-4$, hence Theorem  \ref{thmd3} above is best possible. 
 Related results for rational curve arrangements are stated below in Corollary \ref{corB1} and \ref{corB2}.
\end{rk}

\section{Results involving the Bourbaki ideal and the Jacobian module $N(f)$}

We recall now the construction of the Bourbaki ideal $B(C,r_1)$ associated to a degree $d$ reduced curve $C:f=0$ and to a minimal degree non-zero syzygy $r_1 \in AR(f)$, see \cite{DStJump}.
For any choice of the syzygy $r_1$ with minimal degree $d_1$, we have a morphism of graded $S$-modules
\begin{equation} \label{B1}
S(-d_1)  \xrightarrow{u} AR(f), \  u(h)= h \cdot r_1.
\end{equation}
For any homogeneous syzygy $r=(a,b,c) \in AR(f)_m$, consider the determinant $\Delta(r)=\det M(r)$ of the $3 \times 3$ matrix $M(r)$ which has as first row $x,y,z$, as second row $a_1,b_1,c_1$ and as third row $a,b,c$. Then it turns out that $\Delta(r)$ is divisible by $f$, see \cite{Dmax}, and we define thus a new morphism of graded $S$-modules
\begin{equation} \label{B2}
 AR(f)  \xrightarrow{v}  S(d_1-d+1)   , \  v(r)= \Delta(r)/f,
\end{equation}
and a homogeneous ideal $B(C,r_1) \subset S$ such that $\im v=B(C,r_1)(d_1-d+1)$.
It is known that the ideal $B(C,r_1)$, when $C$ is not a free curve, defines a $0$-dimensional subscheme in $\PP^2$, which is locally a complete intersection, see \cite[Theorem 5.1]{DStJump}. Moreover, one has the following obvious exact sequence
\begin{equation} \label{B3}
 0 \to S(-d_1)  \xrightarrow{u} AR(f) \xrightarrow{v} B(C,r_1)(d_1-d+1) \to 0.
\end{equation}

\begin{prop}
\label{propD}
With the above notation,  $C:f=0$ is a 3-syzygy curve if and only if the corresponding Bourbaki ideal $B(C,r_1)$, for some minimal degree syzygy $r_1$, is a complete intersection.
More precisely, $r_1,r_2,r_3 \in AR(f)$ form a minimal system of generators for $AR(f)$ if and only if $g_2=v(r_2)$ and $g_3=v(r_3)$ generate the ideal $B(C,r_1)$.
\end{prop}
\proof
The proof follows from the exact sequence \eqref{B3}.
\endproof

\begin{rk}
\label{rkcoeff} Let $C:f=0$ be a 3-syzygy curve as above.
Consider the determinant $\Delta(R)=\det M(R)$ of the $3 \times 3$ matrix $M(R)$ which has as first row $x,y,z$, as second row the components of $r_1$, namely $a_1,b_1,c_1$, and as third row the components of $R$, the combination of $r_1,r_2$ and $r_3$ considered in Corollary \ref{corA}. Then $R=0$ implies that $\Delta(R)=0$. On the other hand, expanding this determinant with respect to the third row, we get
$$h_2g_2+h_3g_3=0.$$
Note that $\deg h_2=d_1+d_3-d+1=\deg g_3$ and $\deg h_3=d_1+d_2-d+1=\deg g_2$.
Since the ideal $(g_2,g_3)$ defines a 0-dimensional complete intersection in $\PP^2$, it follows that $g_2$ and $g_3$ have no common factor. We conclude that 
$$h_2=\al g_3 \text{  and } h_3=-\al g_2,$$
 for some $\al \in \C^*$.
\end{rk}

\begin{prop}
\label{propB}
Suppose $C:f=0$ is a 3-syzygy curve of degree $d\geq 3$ with exponents $(d_1,d_2,d_3)$. Then  the total Tjurina number $\tau(C)$ is given by the equivalent formulas
$$\tau(C)=(d-1)^2-d_1(d-d_1-1)-(d_3-d_2+1)-(d_1+d_2-d)(d_1+d_3-d+2)$$
and
$$\tau(C)=(d-1)(d-d_1-1)+d_1^2-[d_1-(d-1-d_2)][d_1-(d-1-d_3)].$$
In particular, the curve $C$ 
satisfies
$$\tau(C)\leq (d-1)^2-d_1(d-d_1-1)-(d_3-d_2+1)$$
and equality holds if and only if $C$ is a plus one generated curve.

\end{prop}

\proof
A direct computation using the formula for $\tau(C)$ given in Proposition \ref{propA} or the equality
\begin{equation}
\label{tau1} 
\deg B(C,r_1)=(d-1)^2-d_1(d-d_1-1)-\tau(C),
\end{equation}
given in \cite[Theorem 5.1]{DStJump}.
Indeed, one clearly has $\deg B(C,r_1)= \deg g_2 \deg g_3=$
$$=(d_1+d_3-d+1)(d_1+d_2-d+1)=(d_3-d_2+1)+(d_1+d_2-d)(d_1+d_3-d+2).$$
\endproof

\begin{cor}
\label{corBour}
Let $C:f=0$ be a reduced curve of degree $d\geq 3$ and let $r_j$ for $j=1,...,m$ be a minimal system of homogeneous generators for $AR(f)$ such that $d_1 \leq \cdots \leq d_m$, where $d_j=\deg r_j$. Assume that $m\geq 3$, i.e. the curve $C$ is not free, and  that $g_k=v(r_k)$ for $k=2,3$ define a complete intersection ideal $B'(C,r_1)$ in $S$. Then
$$\tau(C) \geq (d-1)(d-d_1-1)+d_1^2-[d_1-(d-1-d_2)][d_1-(d-1-d_3)]$$
and equality holds if and only if $C$ is a 3-syzygy curve.
\end{cor}
\proof
The inclusion $B'(C,r_1) \subset B(C,r_1)$ implies that $\deg B'(C,r_1) \geq \deg B(C,r_1)$, and this yields the claimed inequality. On the other hand, the equality $\deg B'(C,r_1) = \deg B(C,r_1)$ implies the equality $ B'(C,r_1) = B(C,r_1)$. Indeed, the equality of degrees implies that 
$\dim B'(C,r_1)_s=\dim B(C,r_1)_s$ for all large $s$. Hence, for any homogeneous $ h \in B(C,r_1)$, there is a large $t$ such that $(x,y,z)^th \in  B(C,r_1)_s= B'(C,r_1)_s$, i.e.
the polynomial $h$ belongs to the saturation $B'(C,r_1)^{sat}$ of $B'(C,r_1)$ with respect to the maximal ideal ${\bf m}=(x,y,z)$. Since $B'(C,r_1)$ is a complete intersection, we have $B'(C,r_1)^{sat}=B'(C,r_1)$, and hence $ B'(C,r_1) = B(C,r_1)$.
\endproof

It was shown by A. du Plessis and C.T.C. Wall that, for any plane curve $C:f=0$ of degree $d$, one has the inequality
\begin{equation}
\label{dPW} 
\tau(C) \geq (d-1)(d-d_1-1),
\end{equation}
where $r=d_1=mdr(f)$, see \cite[Theorem 3.2]{dPW}, and also \cite[Theorem 20]{E} for a new approach. The following result describes the curves for which this lower bound is attained. For explicit examples of such minimal Tjurina curves, see Example \ref{exConj} (i) and Example \ref{ex1.1}.

\begin{thm}
\label{thmNEW}
Let $C:f=0$ be a reduced curve of degree $d\geq 3$ with exponents $d_1 \leq \cdots \leq d_m$ with $m \geq 3$, and let $r_1$ be a non-zero syzygy of minimal degree. Let $d'$ be the smallest integer such that $$d_3 \leq d' \leq \min(d_m,d-1)$$
and the linear system $B(C,r_1)_{d_1+d'-d+1}$ has a 0-dimensional base locus.
Then 
$$
\tau(C) \geq (d-1)(d-d_1-1)+d_1^2-[d_1-(d-1-d_2)][d_1-(d-1-d')],
$$
and equality holds if and only if $C$ is a 3-syzygy curve and then $d'=d_3$.
In particular, one has the following.
\begin{enumerate}

\item The equality $\tau(C)=(d-1)(d-d_1-1)$
holds if and only if $C$ is a 3-syzygy curve and 
$$d_2=d_3=d-1.$$

\item When $C$ is a line arrangement, then
$$\tau(C) \geq (d-1)(d-d_1-1)+2d_1-1.$$
In other words, for line arrangements $C$, the du Plessis-Wall inequality
\eqref{dPW} is strict.

\end{enumerate}
\end{thm}

\proof
Note first that such an integer $d'$ always exists. Indeed, consider the graded $S$-submodule $KR(f) \subset AR(f)$, generated by the Koszul type syzygies $k^x,k^y$ and $k^z$ considered in the proof of Theorem
\ref{thmd3}, and let $B_K(C,r_1)(d_1-d+1)$ denote the image of $KR(f)$ under the morphism $v$ in the exact sequence \eqref{B3}. Then
the ideal $B_K(C,r_1) \subset B(C,r_1)$ defines a 0-dimensional subscheme in $\PP^2$ if and only if the sequence
$\dim \left(B(C,r_1)/B_K(C,r_1)\right)_n$ is bounded. This is indeed the case, since there is an epimorphism
$$\left(AR(f)/KR(f) \right)_n \to \left(B(C,r_1)/B_K(C,r_1)\right)_n$$
induced by $v$, and the sequence $\dim \left(AR(f)/KR(f) \right)_n$ is bounded by $\tau(C)$, see \cite[Corollary 11]{CD} or \cite[Theorem 1]{DBull}.
Since the ideal $B_K(C,r_1)$, which is generated by elements of degree
${d_1}$, defines a $0$-dimensional subscheme in $\PP^2$, it follows that $B(C,r_1)_{d_1}$ has a 0-dimensional base locus.
We choose $d'$ the smallest positive integer $\geq d_3$ such that
$B(C,r_1)_{d_1+d'-d+1}$ has also a 0-dimensional base locus. Clearly one has $d'\leq d_m$.
Then choose $r \in AR(f)_{d'}$ such that $g=v(r)$ is a generic element in $B(C,r_1)_{d_1+d'-d+1}$. It follows that $g_2$ and $g$ defines a 0-dimensional complete intersection ideal $B''(C,r_1)\subset B(C,r_1)$. The same argument as above in the proof of Corollary \ref{corBour} completes the proof. 
Note that both factors $d_1-(d-1-d_2)$ and $d_1-(d-1-d')$ are in the interval $[1,d_1]$, and hence their product is $\leq d_1^2$.
For the last claim, recall that one has 
 $d_m \leq d-2,$
see \cite[Corollary 3.5]{Sch0}, and hence our integer $d'$ defined above is also $\leq d-2$. This implies that for a line arrangement $C$ one has
$$\tau(C) \geq (d-1)(d-d_1-1)+d_1^2-[d_1-(d-1-d_2)][d_1-(d-1-d')]  \geq $$
$$\geq (d-1)(d-d_1-1)+d_1^2-(d_1-1)^2=(d-1)(d-d_1-1)+2d_1-1.$$
\endproof

\begin{ex}
\label{exconj1} 
With the above notation, the ideal $B'(C,r_1)=(g_2,g_3)$ is a complete intersection for many
curves, see Example \ref{exConj}. However, there are curves for which the ideal $B'(C,r_1)$ is a not a complete intersection.
For instance, the curve $C:f=0$, with 
$$f= (x^2+y^2)^3-4x^2y^2z^2=0,$$
 has 4 syzygies of degrees 3,4,4,4 as a minimal set of generators for $AR(f)$, namely
$$r_1=(-xy^2,x^2y, -x^2z+y^2z), \  r_2=(-3y^4, 3xy^3-4xyz^2,-3x^3z-9xy^2z+4xz^3),$$
$$r_3=(-3x^3y+4xyz^2,3x^4,9x^2yz+3y^3z-4yz^3)$$ 
{ and } 
$$r_4=(4xy^2z,0,3x^4+6x^2y^2+3y^4-4y^2z^2).$$
It follows that $g_2=v(r_2)=-3yz$, $g_3=v(r_3)=3xz$ and $g_4=3xy$.
Hence
in this case the polynomials  $g_2$ and $g_3$ are divisible by $z$. However, if one replces $r_3$ by $r_3'=ar_3+br_4$ for $a,b \in \C^*$, then $g_3'=v(r_3')=3x(az+by)$, and hence $g_2$ and $g_3'$ have no common factor, which implies that $d'=4=d_3$ in this case.
Note that for this curve $C$ one has $ \tau(C)=16$ and $\nu(C)=3$, and the minimal resolution for $M(f)$ has the form
$$
0 \to S(-10)^2 \to S(-8) \oplus S(-9)^3 \to S(-5)^3 \to S.$$
The curve $C$ has 3 singularities: two simple cusps $A_2$ located at $(1:\pm i:0)$, and a singularity with Milnor number 13 and two branches located at $(0:0:1)$. It follows that $C$ has genus 1.

If we consider the rational  curve $C': f'= (x^2+y^2)^6-3(x^{11}+y^{11})z=0$, then a computation by SINGULAR \cite{Sing} shows that $AR(f')$ has a minimal set of generators consisting of 4 syzygies, of degrees 6,7,8 and 10 respectively and too complicated to write down, and the corresponding polynomials  $g_2$ and $g_3$ are divisible by $z$ again, but $g_4$ is an irreducible quintic. Hence in this case we have $d'=8>d_3$.
These two curves $C$ and $C'$  show that curves with $d_1+d_2=d+1$ can be rather complicated, unlike the plus-one curves corresponding to $d_1+d_2=d$.

\end{ex}

 It was shown in \cite[Corollary 4.3]{DPop} that the graded $S$-module  $N(f)$ satisfies a Lefschetz type property with respect to multiplication by generic linear forms. This implies in particular the inequalities
\begin{equation}
\label{in} 
0 \leq n(f)_0 \leq n(f)_1 \leq ...\leq n(f)_{[T/2]} \geq n(f)_{[T/2]+1} \geq ...\geq n(f)_T \geq 0,
\end{equation}
where $T=3d-6=$ end$(M(f))$ in the notation from \cite{HS}. We have the following result, which generalizes the fact that $\nu(C)=1$  and $ct(f)+st(f)=T+2$ for a nearly free curve $C:f=0$.
\begin{prop}
\label{propC} Suppose $C:f=0$ is a 3-syzygy curve of degree $d\geq 3$ with exponents $(d_1,d_2,d_3)$ such that $d_1 \leq d/2$. Then
$$\nu(C)=d_3-d_2+1+(d_1+d_2-d)(d_1+d_3-d+2)=[d_1-(d-1-d_2)][d_1-(d-1-d_3)].$$
In particular, when $C:f=0$ is a plus-one generated curve of degree $d\geq 3$ with exponents $(d_1,d_2,d_3)$, then
$$\nu(C)=d_3-d_2+1 \leq d_1 \text{ and } ct(f)+st(f)=T+ \nu(C)+1.$$
\end{prop}
\proof
Use \cite[Theorem 1.2]{Drcc}, with the remark that the equality in (1) there holds for 
$r=d_1=mdr(f) \leq d/2$, as well as
Proposition \ref{propA}, Proposition \ref{propB} and Theorem \ref{thmd3} above.
\endproof

\begin{ex}
\label{exTS}
For a Thom-Sebastiani curve $C:g(x,y)+z^d=0$, where $g$ has distinct $m$ factors as in Example \ref{ex1.1} below, we get $ct(f)=d+m-3$, $st(f)=2d+m-5$ and hence $ct(f)+st(f)=T+2m-2$.
On the other hand one has
$\nu(C)=m(m-2)+1$, and hence a linear relation involving $ct(f)$, $st(f)$ and $\nu(C)$ does not seem to exist for 3-syzygy curves.
\end{ex}

We also  consider the invariant 
$$\sigma(C)=\min \{j   : n(f)_j \ne 0\}= \indeg (N(f)),$$
in the notation from \cite{HS}.
The self duality of the graded $S$-module $N(f)$, see \cite{HS,Se, SW}, implies that 
\begin{equation}
\label{dual} 
 n(f)_{k} = n(f)_{T-k},
\end{equation}
for any integer $k$, hence
$n(f)_s \ne 0$ exactly for $s=\sigma(C),..., T-\sigma(C)$. Hassanzadeh and Simis result in \cite[Proposition 1.3]{HS} with our slight improvement in \cite[Remark 4.2]{DStSat}, implies the following.
\begin{thm}
\label{thmHS}
Let $C:f=0$ be a 3-syzygy curve with exponents $(d_1,d_2,d_3)$ and set  $e=d_1+d_2+d_3$ as above. Then the following hold.

\medskip

\noindent (i) The minimal free resolution of $N(f)$ as a graded $S$-module has the form
$$0 \to  S(-e)\to \oplus_{j=1} ^3S(1-d-d_j)  \to \oplus_{j=1} ^3S(d_j+2-2d) \to S(e+3-3d),$$
where the leftmost map is the same as in the resolution \eqref{res2}. In particular,
$$\sigma(C)= \indeg (N(f))=3(d-1)-(d_1+d_2+d_3).$$

\medskip

\noindent (ii) If $d_1 \geq 2$, then the resolution of $S/I_f$ is given by
$$0  \to \oplus_{j=1} ^3S(d_j+2-2d) \to S^3(1-d)\oplus S(e+3-3d) \to S.$$

\end{thm}

\begin{cor}
\label{corHS}
Suppose $C:f=0$ is a plus-one generated curve of degree $d\geq 3$ with exponents $(d_1,d_2,d_3)$, which is not nearly free, i.e. $d_2<d_3$. Set $k_j=2d-d_j-3$ for $j=1,...,3$. 
Then one has the following minimal free resolution of $N(f)$ as a graded $S$-module.
$$0 \to S(-d-d_3) \to S(-d-d_3+1) \oplus S(-k_1-2) \oplus S(-k_2-2) \to $$
$$\to S(-k_1-1) \oplus S(-k_2-1) \oplus S(-k_3-1) \to S(-k_3).$$
In particular, $\sigma(C)=k_3$ and the Hilbert function of $N(f)$ is given by the following formulas.
\begin{enumerate}

\item $n(f)_k=0$ for $k <k_3$, and $k_3<T/2$.

\item $n(f)_k=k-k_3+1$ for $k_3 \leq k \leq k_2$, and $k_2=d+d_1-3 \leq T/2$.

\item $n(f)_k=d_3-d_2+1=\nu(C)$ for $k_2 \leq k \leq T/2$.

\end{enumerate}

\end{cor}
Note that the dimensions $n(f)_k$ for $k>T/2$ can be determined from this result using the equality \eqref{dual}. When $C:f=0$ is a nearly free curve, the numbers $n(f)_k$ can be given by the above formulas as well, but we have to take $k_3=k_2 \leq T/2$, in other words we have only the cases (1) and (3).

\begin{thm}
\label{thmPO2}
Let $C:f=0$ be a reduced, non-free plane curve of degree $d$ such that the Hilbert function of the Jacobian module $N(f)$ has the form described in Corollary \ref{corHS}, namely suppose that there are integers $0 <k_3=\sigma(C) \leq k_2 \leq T/2$ such that 
\begin{enumerate}

\item $n(f)_k=0$ for $k <k_3$.

\item $n(f)_k=k-k_3+1$ for $k_3 \leq k \leq k_2$.

\item $n(f)_k=k_2-k_3+1=\nu(C)$ for $k_2 \leq k \leq T/2$.
\end{enumerate}
Then $C$ is a plus-one generated curve. Moreover, a reduced plane curve $C$ of degree $d$ such that $$\nu(C)=2$$
 is one of the following.
\begin{enumerate}

\item[(i)] a plus-one generated curve with exponents $(d_1,d-d_1,d-d_1+1)$ and total Tjurina number
$$\tau(C)= (d-1)^2-d_1(d-d_1-1)-2;$$

\item[(ii)] a 4-syzygy curve with exponents $d_1=d_2=d_3=d_4$, degree $d=2d_1-1$ and
total Tjurina number 
$$\tau(C)=3d_1^2-6d_1+1=(d-1)^2-d_1(d-d_1-1)-3.$$ 
When $d=3$, then $C$ is a nodal cubic.
\end{enumerate} 
 \end{thm}

\proof
When $C$ is a plus-one generated curve, then the formulas for the exponents and for the 
total Tjurina number follow from Proposition \ref{propC} and Proposition \ref{propA} (4).
We assume now that $C$ is not a plus-one generated curve, which in view of Theorem \ref{thmPO1} means exactly that $d_1+d_2>d$.
We prove that this inequality implies that the Hilbert function of the Jacobian module $N(f)$ has not the form described in Corollary \ref{corHS} and that $\nu(C) \geq 3$, unless $C$ is as in (ii) above. Using the notation from the resolution \eqref{res2A}, we assume first that $m=3$, i.e. $C$ is a 3-syzygy curve.
Consider the resolution of $N(f)$ given in Theorem \ref{thmHS} (i) and take the graded pieces of degree $s=\sigma(C)+1=3(d-1)-(d_1+d_2+d_3)+1$. Note that the inequality $d_1+d_2>d$
implies that 
$$s+d_j+2-2d \leq s+d_3+2-2d=d-d_1-d_2<0.$$
It follows that $n(f)_s=3$, and hence our claims are proved in this case. 

Assume now that $m>3$
and consider the resolution of $N(f)$ obtained from the resolution \eqref{res2A} by applying
 \cite[Proposition 1.3]{HS}, namely the resolution
 \begin{equation}
\label{res20A}
0 \to \oplus_{i=1} ^{m-2}S(-e_i) \to \oplus_{j=1} ^mS(1-d-d_j)\to  \oplus_{j=1} ^mS(d_j-2(d-1))   \to  \oplus_{i=1} ^{m-2}S(e_i-3(d-1)).
\end{equation}
Using this resolution, denoted shortly by $0 \to F_3 \to F_2 \to F_1 \to F_0$,  it follows that
\begin{equation}
\label{sigma}
\sigma(C)=\indeg (N(f))=3(d-1)-e_{m-2}=2(d-1)-d_m-\epsilon_{m-2},
\end{equation}
with the notation from \eqref{res2B}.
Note that $\min(2(d-1)-d_j)=2(d-1)-d_m >\sigma(C)$, and hence using the component of degree
$\sigma(C)$ of the above resolution, we conclude that $F_{1,\sigma(C)}=0$ and $n(f)_{\sigma(C)}=t$, where $t$ is the cardinal of the set $\I=\{i \ : \ e_i=e_{m-2}\}$.
If $t>1$, then there are two cases to discuss.

\medskip

\noindent{ \bf Case 1.}  Suppose that one has $1 <t<m-2$.
Then consider the homogeneous component of degree $s=\sigma(C)+1$ of this resolution. It follows that $\dim F_{0,s} \geq 3t$, with equality when there is no $i\notin \I$ with $e_i=e_{m-2}-1$, and $\dim F_{1,s} \leq t$, with equality when $\epsilon_i=1$ for all $i \in \I$.
It follows that $n(f)_s \geq 2t \geq 4$ for $t\geq 2$. 

\medskip

\noindent{ \bf Case 2.}  Suppose that one has $1 <t=m-2$. Note that in this case 
$\dim F_{1,s} \leq m$, with equality when $d_1=...=d_m$ and $\epsilon_i=1$ for all $i$, and hence
$n(f)_s \geq 3(m-2)-m =2m-6\geq 2$. Therefore, one has $n(f)_s>2$ except the case when
$d_1=...=d_m$,  $\epsilon_i=1$ for all $i$, and $m=4$. In this case the equality \eqref{res2C} implies that
\begin{equation}
\label{E1}
2d_1=d+1.
\end{equation}
Since $d \geq 3$ it follows that $d_1 \geq 2$. For $d_1=2$ we get $d=3$, and the classification of the cubic plane curves shows that $C$ must be a nodal cubic.
We  consider now the case $d_1 \geq 3$. To do this we use \cite[Proposition 1.3]{HS} completed by our \cite[Remark 4.2]{DStSat}, and get the following resolution for $S/I_f$
\begin{equation}
\label{E2}
0 \to S^m(2(1-d)+d_1) \to S^3(1-d) \oplus S^{m-2}(d_1+3-2d) \to S.
\end{equation}
Using the formula for $\tau(C)$ in terms of such a resolution, see for instance \cite[Lemma 2.2 (4)]{DStSat}, we get the following
\begin{equation}
\label{E3}
2\tau(C)=4(2(d-1)-d_1)^2-3(d-1)^2-2(2d-3-d_1)^2.
\end{equation}
Using the formula \eqref{E1}, this gives
$$\tau(C)=3d_1^2-6d_1+1.$$

\medskip

We  suppose from now on that $t=1$ and hence
$e_{m-3}<e_{m-2}$. To fix the ideas, suppose that $e_{m-3}=e_{m-2}-a$, for an integer $a>0$.
 It follows that, for $s=\sigma(C)+1$ as above, we have $\dim F_{0,s} \geq 3$, with equality when $a>1$. 
 Hence, in order to get $n(f)_s \leq 2$, we need
 $F_{1,s}\ne 0$, which implies $\epsilon_{m-2}=1$.
 Then the equality $e_{m-3}=e_{m-2}-a$ implies $d_{m-1}\leq d_{m}-a$. The homogeneous component of degree $s'=\sigma(C)+a$ of the resolution is
 $$0 \to F_{3,s'} \to F_{2,s'} \to F_{1,s'} \to F_{0,s'}$$
 such that the following hold.
 \begin{enumerate}

\item[(i)] $\dim F_{0,s'} \geq \dim S_a+\dim S_0$, where $S_a$ corresponds to $S(e_{m-2}-3(d-1))_{s'}$ and
$S_0$ to $S(e_{m-3}-3(d-1))_{s'}$.

\item[(ii)] $\dim F_{1,s'}=\dim S_{a-1}$, where $S_{a-1}$ corresponds to $S(d_{m}-2(d-1))_{s'}$, and using
the inequality $d_{m-1}-2(d-1)+s' <0$.

\end{enumerate}
This implies that 
$$n(f)_{\sigma(C)+a} \geq \dim S_a+\dim S_0-\dim S_{a-1}=a+2 \geq 3,$$
and hence the Hilbert function of the Jacobian module $N(f)$ has not the form described in Corollary \ref{corHS} and also $\nu(C)\geq 3$.
\endproof

\begin{ex}
\label{exNU=2} The two types of curves having $\nu(C)=2$ in Theorem \ref{thmPO2} (ii) correspond to the two cases in \cite[Corollary 23 (iii)]{E0}. To relate the two results, one has to use the rank two vector bundle $E_C$ which corresponds to the graded $S$-module $AR(f)$, which is the main object of study in \cite{AD,DStJump}.

Here are some examples of 4-syzygy curves $C$ with $\nu(C)=2$ as described in Theorem \ref{thmPO2} (ii).
The generic line arrangement
$$\A_5: f_5=xyz(x-2y-3z)(x+y+z)=0,$$
the line arrangement
$$\A_7: f_7= xyz(x - z)(y - z)(x + y)(x+y+z)=0 $$
and the line arrangement
$$\A_9: f_9= xyz(x - z)(y - z)(x + y)(x+y+z)(x-y)(x-y-z)=0  $$
are all of this type. It is interesting to note that the line arrangement
$$\B_7: g_7=xyz(x+y+z)(x+z)(y+z)(x-y+2z)=0$$
is a plus-one generated line arrangement with exponents $(d_1,d_2,d_3)=(3,4,5)$, and the arrangements $\A_7$ and $\B_7$ have both 9 double points and 4 triple points as multiple points. In particular they have $\tau(\A_7)=\tau(\B_7)=25$ and $\nu(\A_7)=\nu(\B_7)=2$. However, the line arrangements 
$\A_7$ and $\B_7$ are combinatorially distinct, since any line in $\A_7$ contains at least a triple point, but the line $x-y+2z=0$ in $\B_7$ contains no triple point.

Beyond the case of line arrangements, note that the reducible curve 
$$C'_5:f'_5=(x + y)(x^4 - x^3y + x^2y^2 - xy^3 + y^4 + x^3z - x^2yz + xy^2z)=0$$
and the irreducible curve
$$C''_5:f''= z(x^2+y^2+xy)^2+x^5+y^5=0$$
are both 4-syzygy curves as in Theorem \ref{thmPO2} (ii) with $d_1=3$. A curve similar to $C''_5$ is studied in \cite[Example 6.3]{DStJump}. The irreducible curve
$$C'''_9:f'''_9= z(x^2+y^2+xy)^4+x^9+y^9=0$$
is a 4-sygygy curve as in Theorem \ref{thmPO2} (ii) with $d_1=5$.
The curves $C''_5$ and $C'''_9$ are both rational, nearly cuspidal curves: the only singularity is located at $(0:0:1)$ and has two branches.
Note also that there are curves $C$ of even degree with $\nu(C)=3$ which are 4-syzygy curves, as shown by the first curve in Example \ref{exconj1}. 

\end{ex}

\begin{cor}
\label{corHS2}
Let $C:f=0$ be a reduced plane curve of degree $d\geq 3$. Then the fact that $C$ is a plus-one generated curve, and in the affirmative case, the corresponding exponents $(d_1,d_2,d_3)$, are determined by the degree $d$ and the Hilbert function of the Jacobian module $N(f)$.
\end{cor}

\section{Examples of 3-syzygy curves and plus one generated curves} 

First we look at some low degree, irreducible 3-syzygy curves and plus one generated curves.

\begin{ex}
\label{exlowdegree}
Here are some classical plane curves.

\medskip

\noindent (i) The curve  $C:f =(y^2-xz)^2+y^2z^2+z^4=0$ has a unique singularity of type $A_5$ located at $p=(1:0:0)$, and it is a plus one generated curve with exponents $(d_1,d_2,d_3)=(2, 2, 3)$.
One has $\delta(C,p)=3$, hence $C$ is a rational curve with a point $p$ having $r(C,p)=2$ branches.

\medskip

\noindent (ii) The curve  $C:f =(x^2+y^2)^2-4xy^2z=0$, called a {\it double folium},  has a unique singularity of type $D_5$ located at $p=(0:0:1)$, and  $C$ is a plus-one generated curve with exponents $(d_1,d_2,d_3)=(2, 2, 3)$.
One has again $\delta(C,p)=3$, hence $C$ is again a rational curve with a point $p$ having $r(C,p)=2$ branches.

\medskip

\noindent (iii) The curve  $C:f =(x^2+y^2-2xz)^2-(x^2+y^2)z^2=0$, called a {\it lima\c con},  has 3 singularities, one of type  $A_1$ located at $p_1=(0:0:1)$, and two of type $A_2$, located at $(1:\pm i:0)$. The curve  $C$ is a plus one generated curve with exponents $(d_1,d_2,d_3)=(2, 2, 3)$.
One has again $\delta(C,p)=3$, hence $C$ is again a rational curve with a point $p_1$ having $r(C,p_1)=2$ branches.

\medskip

\noindent (iv) The curve  $C:f =x^5-y^2z^3-xz^4=0$, called a {\it Bolza curve},  has a unique singularity of type $E_8$ located at $p=(0:1:0)$, and  $C$ is a 3-syzygy curve with exponents $(d_1,d_2,d_3)=(2, 4, 4)$.
One has $\delta(C,p)=4$, hence $C$ is  cuspidal curve of genus $g=2$.

\medskip

\noindent (v) The curve  $C:f =x^6+y^6-x^2z^4=0$, called a {\it Butterfly curve},  has a unique singularity of type $A_5$ located at $p=(0:0:1)$, and  $C$ is a 3-syzygy curve with exponents $(d_1,d_2,d_3)=(4, 5, 5)$.
One has $\delta(C,p)=3$, hence $C$ is  a curve of genus $g=7$,
with a point $p$ having $r(C,p)=2$ branches.

\end{ex}

\begin{ex}
\label{exNUlarge}
Consider the rational curve 
$$f=x^{2d'+1}+(x^{d'}+y^{d'})^{2}z=0.$$
It can be shown that $C$ has the exponents $(d',d'+1,2d')$ for any $d'\geq 3$, hence we have a plus-one generated curve
with $\nu(C)=d_3-d_2+1=d'$ arbitrarily large.
Indeed, the corresponding generators for $AR(f)$ are the following.
$$r_1=(0,x^{d'}+y^{d'},-2d'y^{d'-1}z),$$
$$r_2=(-2d'(x^{d'}+y^{d'})z, dy^{d'+1}, 2d'dx^{d'}z-2d'dy^{d'}z+4(d')^2x^{d'-1}z^2),$$
and $r_3=(f_z,0,-f_x)$. Proposition \ref{propC} says that for a plus-one generated curve $C$ we have $\nu(C) \leq d_1$, and here we have a sequence of examples of plus-one generated curves for which this inequality is in fact an equality.
\end{ex}

Next we present two infinite series of 3-syzygy curves, each having 2 irreducible components: a line and a smooth curve.
\begin{ex}
\label{exConj}
Let $C':f'=0$ be a smooth curve of degree $d-1\geq 2 $ and consider the curve $C:f(x,y,z)=xf'(x,y,z)=0$, namely $C$ is the union of the smooth curve $C'$ with the line $L:x=0$.

\medskip

\noindent (i) Consider first the case $f'=x^{d-1}+y^{d-1}+z^{d-1}$. Then $L$ is transverse to $C'$, $C$ has $d-1$ simple nodes as singularities and hence $\tau(C)=d-1$. The module is generated by 3 syzygies, namely
$$r_1=(0,f'_z,-f'_y), \ r_2=(f_y,-f_x,0) \text{ and } r_3=(f_z,0,-f_x)$$
of degrees $d_1=d-2$ and $d_2=d_3=d-1$. Hence the du Plessis-Wall inequality in \eqref{dPW} is an equality for this curve.
The ideal $B'(C,r_1)$ coincides with the ideal
$(f'_y,f'_z)$ and hence it is a complete intersection. In conclusion $C$ is a 3-syzygy curve for any $d\geq 3$, and a plus one generated curve for $d=3$.

\medskip

\noindent (ii)  Consider now the case $f'=x^{d-1}+xy^{d-2}+z^{d-1}$. Then $L$ meets $C'$ only at the point $p=(0:1:0)$, where $C$ has an $A_{2d-3}$ singularity, and hence $\tau(C)=2d-3$. The module is generated by 3 syzygies, namely
$$r_1=(0,f'_z,-f'_y), \ r_2=(-xy^{d-3}, \frac{dx^{d-2}+2y^{d-2}}{d-2},\frac{y^{d-3}z}{d-1}) \text{ and } r_3=(f_z,0,-f_x)$$
of degrees $d_1=d_2=d-2$ and $d_3=d-1$. In fact, for $d=3$, the last syzygy is a consequence of the first two, and we get a free curve. For $d=4$ we get a plus one generated curve,
and for any $d>4$ we get a 3-syzygy curve.

\medskip

\noindent (iii) Finally we consider the case $f'=x^4+xy^3+xz^3+y^2z^2$. Then L is a simple bitangent, and hence $C$ has two singularities $A_3$, which give $\tau(C)=6$.
A computation using Singular shows that $AR(f)$ has 4 minimal generators of degrees
$d_1=3$, $d_2=d_3=d_4=4$, namely
$$r_1=(0,-2y^2z-3xz^2,3xy^2+2yz^2),$$
$$r_2=(-6xy^3+6xz^3,10x^3y+4y^4+yz^3,-10x^3z-y^3z-4z^4),$$
$$r_3=(-3x^2y^2-2xyz^2,5x^4+2xy^3+y^2z^2+2xz^3,0),$$
and
$$r_4=(-6xy^2z-9x^2z^2,4y^3z+6xyz^2,15x^4-y^2z^2+6xz^3).$$
Hence this curve $C$ is not a 3-syzygy curve.
The corresponding ideal $B'(C,r_1)$ is the complete intersection 
$(g_2,g_3)=(y^3-z^3,3xy+2z^2)$ and the inequality in Corollary \ref{corBour} is in this case
$\tau(C)=6 >4.$ 
\end{ex}

We present next an infinite series of plus one generated irreducible curves.
\begin{ex}
\label{ex2}
Consider the irreducible curve $C:f=0$ with  $d=2k+1$ odd, $k \geq 2$, given by
$$f=x^d+(x^2+y^2)^{k}z=0.$$
Then the exponents are $(2,d-2,d-1)$, hence we have a  plus-one generated curve. The generating syzygies can be described as follows. Note that 
$$f_x=(2k+1)x^{2k}+2k(x^2+y^2)^{k-1}xz, \ \ f_y=2k(x^2+y^2)^{k-1}yz \text{ and } f_z=(x^2+y^2)^{k}.$$
The first syzygy $r_1$ of degree $d_1=2$ is easy to find, namely
$$(x^2+y^2)f_y-2kyzf_z=0.$$
The last syzygy $r_3$ of degree $d-1$ can be taken to be the Koszul syzygy
$$(x^2+y^2)^{k}f_x-((2k+1)x^{2k}+2k(x^2+y^2)^{k-1}xz)f_z=0.$$
The second syzygy $r_2$ of degree $d-2=2k-1$ is more complicated to describe. We consider
first the case $k\geq 3$ odd, and set
$$A=af_x+bf_y$$
where $a=2kz(x^2+y^2)^{k-1}$ and $b=(2k+1)y^{2k-1}$.
Then we have
$$A=2k(2k+1)z(x^{2k}+y^{2k})(x^2+y^2)^{k-1}+4k^2xz^2(x^2+y^2)^{2k-2}.$$
Since $k$ is odd, we have 
$$x^{2k}+y^{2k}=(x^2+y^2)(x^{2k-2}-x^{2k-4}y^2+ ...-x^2y^{2k-4}+y^{2k-2}),$$
and hence one can write $A=Bf_z$, where
$$B=2k(2k+1)z(x^{2k-2}-x^{2k-4}y^2+ ...-x^2y^{2k-4}+y^{2k-2})+4k^2xz^2(x^2+y^2)^{k-2}.$$
With this notation, the syzygy $r_2$ is given by
$$af_x+bf_y-Bf_z=0.$$
When $k$ is even, one can do essentially the same trick, starting with $A=af_x-bf_y$.
To prove that these 3 syzygies generate $AR(f)$, one can use Proposition \ref{propD} above.
The formula in Proposition \ref{propA} implies
$$\tau(C)=(d-1)(2d-1)-2(d-2)-2(d-1)-(d-1)(d-2)=d^2-4d+5.$$
This curve is clearly rational, and it has only one singular point, at $p=(0:0:1)$, with Milnor number $\mu(C,p)= d^2-3d+1$ and delta invariant $\delta (C,p)=(d-1)(d-2)/2$. It follows that the germ $(C,p)$ has 
$$r(C,p)=2\delta (C,p)-\mu(C,p)+1=2$$
branches, each of multiplicity $k$.
Using \cite[Theorem 1.2]{Drcc}, it follows that
$$\nu(C)=(d-1)^2-2(d-3)-\tau(C)=2.$$
More precisely, one has $\sigma(C)=3(d-1)-(2d-1)=d-2$ and
$n(f)_{d-1}=2$ using the resolution in Theorem \ref{thmHS} (i). It follows that
$n(f)_k=0$ for $k \notin [d-2,4d-4]$, $n(f)_{d-2}=n(f)_{2d-4}=1$ and $n(f)_{k}=2$ for
$d-1\leq k \leq 2d-5$.
\end{ex}

\begin{ex}
\label{ex1.1} Let $C: f=0$ be a Thom-Sebastiani plane curve, i.e. a curve such that $f(x,y,z)=g(x,y)+z^d$, where
$g$ is a homogeneous polynomial of degree $d$ in $R=\C[x,y]$. Assume that
$$g=\ell_1^{k_1} \cdots \ell_m^{k_m},$$
where the linear forms $\ell_j \in R$ are distinct. Then it is easy to see that the minimal resolution of the Milnor algebra $M(g)=R/(g_x,g_y)$ is given by
$$0 \to R(2-m-d) \to R^2(1-d) \to R.$$
Similarly, $M(z^d)=\C[z]/(z^{d-1})$ has a resolution
$$0 \to \C[z](1-d) \to \C[z].$$
Since $M(f)=M(g) \otimes M(z^d)$, we get the following minimal resolution for the Milnor algebra $M(f)$:
$$0 \to S(3-m-2d) \to S(2-m-d) \oplus S^2(2-2d) \to S^3(1-d) \to S.$$
It follows that $C$ is a 3-syzygy curve with exponents $d_1=m-1$ and $d_2=d_3=d-1$.
Using Corollary \ref{corBour} it follows that
$$\tau(C)=(d-1)^2-d_1(d-d_1-1)-(d_1-1)(d_1+1)=(d-1)(d-d_1-1).$$

Hence the Thom Sebastiani plane curves $C$ discuss in this example realize the minimal value for the global Tjurina number for any pair of given integers $(d,d_1)$. When $d_1=d-2$, this minimal value is also attained when $C$ is the union of a smooth curve of degree $d-1$ and a generic secant line, see \cite[Lemma 21]{E}, which was a motivation for this example. A particular case of such curves occurs in our Example \ref{exConj} (i) above.
\end{ex}

\section{Various bounds on the exponents and the relation to nearly cuspidal rational curves} 
Let $U=\PP^2 \setminus C$ be the complement of $C$ and assume that the curve $C$ has $k$ irreducible components $C_i$, for $i=1,...,k$. Let $g(C_i)$ denote the genus of the normalization of $C_i$. 
Then a key result of U. Walther in \cite[Theorem 4.3]{Wa},  combined with \cite[Theorem 2.7]{Abd}
yields the following inequality.
\begin{thm}
\label{UWthm} For any reduced plane curve $C:f=0$, one has
$$\dim N(f)_{2d-2-j} \leq \dim Gr_F^1H^2(F_f, \C)_{\lambda},$$
for $j=1,2,...,d$ and  $\lambda= \exp(2\pi i (d+1-j)/d)=\exp(-2\pi i (j-1)/d)$. 
In particular, this inequality gives  for  $j=1$, the following
$$n(f)_{d-3}=n(f)_{2d-3} \leq \dim Gr_F^1H^2(F_f, \C)_1= \dim Gr_F^1H^2(U, \C)=\sum_{k=1}^rg(C_i).$$

\end{thm}
\noindent Here $F_f$ is the Milnor fiber of $f$, that is the smooth affine surface in $\C^3$ given by $f(x,y,z)=1$, $F$ denotes the Hodge filtration on the cohomology group $H^2(F_f, \C)$, and
$H^2(F_f, \C)_{\lambda}$ is the Milnor monodromy eigenspace corresponding to the eigenvalue $\lambda$. It is known that one has $H^k(F_f, \C)_1= H^k(U, \C)$ for any $k$, see for instance \cite[Chapter 7]{DHA} for more details on all these objects.
\begin{cor}
\label{corB1}
Let $C:f=0$ be a reduced plane curve of degree $d\geq 3$ with exponents  $d_1 \leq ...\leq d_m$. If all the irreducible components $C_i$ of $C$ are rational, then $d_m\leq d-1$.
\end{cor}
Note that this inequality is sharp even in the class of plus-one irreducible rational curves, as shown by Example \ref{ex2}.
\proof
We have seen in the proof of Theorem \ref{thmPO2} that one has the formula \eqref{sigma}, namely
$$\sigma(C)=2(d-1)-d_m-\epsilon_{m-2}.$$
Theorem \ref{UWthm} implies that $d-3 < \sigma(C)$, which yields our claim since $\epsilon_{m-2}\geq 1$.
\endproof

\begin{cor}
\label{corB2}
Let $C:f=0$ be a 3-syzygy curve of degree $d\geq 3$ with exponents $(d_1,d_2,d_3)$. 
If all the irreducible components $C_i$ of $C$ are rational and $C$ is not a plus-one generated curve, then $d_3\leq d-2$.

\end{cor}

\proof For a 3-syzygy curve we have two formulas for $\sigma(C)$, namely the one in Theorem
\ref{thmHS} and the one in formula \eqref{sigma}.
The equality
$$3(d-1)-d_1-d_2-d_3=2(d-1)-d_3-\epsilon_1,$$
implies $\epsilon_1=d_1+d_2-(d-1)$. If $C$ is not a plus-one curve, then $d_1+d_2>d$ and hence $\epsilon_1\geq 2$. This yields the claim as in the proof above.

\endproof

\begin{thm}
\label{thm2}
Let $C:f=0$ be a rational, nearly cuspidal curve of degree $d$ in $\PP^2$. Assume that one of the following holds. 

\begin{enumerate}

\item[(1)] $d$ is even;

\item[(2)] $d$ is odd and for any singularity $x$ of $C$, the order of any eigenvalue $\lambda_x$ of the local monodromy operator $h_x$ is not $d$.

\item[(3)] $d=p^k$ is a power of a prime number $p$.

\end{enumerate}

Then $\nu(C) \leq 2$. In particular, $C$ is  a free curve when $\nu(C)=0$,  a nearly free curve when $\nu(C)=1$, and a plus-one generated curve when $\nu(C)=2$ and $d$ is even.

\end{thm}

\proof
Note that a rational nearly cuspidal curve $C$ is homotopically equivalent to a sphere $S^2$ with an interval $I=[0,1]$ attached at its extremities. It follows that the Euler characteristic
 of $C$ is given by 
 $$E(C)=E(S^2)+E(I)-E(2 \text{ points })=2+1-2=1.$$
 Then $E(U)=E(\PP^2)-E(C)=3-1=2$ and exactly the same proof as in \cite[Theorem 3.1]{DStRIMS} shows that $\nu(C)\leq 2$ by using Theorem \ref{UWthm} above.
 Then we conclude using our Theorem \ref{thmPO2}.
 \endproof
The fact that $C$ is rational is essential for Theorem \ref{thm2}, as the curve in Remark \ref{rkconj1} shows. Finally we consider in more details the case of rational, nearly cuspidal curve of degree  $d$ odd. 

\begin{thm}
\label{thm3}
Let $C:f=0$ be a rational, nearly cuspidal curve of odd degree $d=2d'+1$. Then $mdr(f)=d_1 \leq d'+1$ and the following hold.

\begin{enumerate}

\item[(1)] If $mdr(f)=d_1 = d'$, then $C$ is  a free curve,  a nearly free curve, or a plus-one generated curve with $\nu(C)=2$;

\item[(2)] If $mdr(f)=d_1 = d'+1$, then $C$ is a 4-syzygy curve as described in Theorem \ref{thmPO2} {\rm (ii)}.

\end{enumerate}

In particular, when  $mdr(f)=d_1\in [d',d'+1]$, then  $\nu(C) \leq 2$.

\end{thm}

\proof
The proof given for \cite[Theorem 1.1]{DStMos} applies with no modification.

\endproof

Using this result, we can construct the following infinite series of nearly cuspidal rational curves
which are 4-syzygy curves as in (2) above.
\begin{ex}
\label{exRNC}  Let $d =2r-1\geq 5$ be an odd integer and set 
$$C_d:f_d=x^{r-1}y^{r-1}z+x^d+y^d=0.$$
 A minimal degree syzygy for $f_d$ is given by
$$\rho_1=(-xy^{r-1},0, (r-1)y^{r-1}z+(2r-1)x^r),$$ 
and hence $d_1=mdr(f_d)=r$. 
The curve $C_d$ is clearly rational, since we can express $z$ as a rational function of $x$ and $y$. The Milnor number $\mu(C_d,p)$ can be easily computed, since the singularity $(C_d,p)$ is  Newton nondegenerate and commode, see \cite{K}. It follows that
$$\mu(C_d,p)=4r^2-10r+5.$$
Since $C_d$ is rational, we have for the $\delta$-invariant the following equality
$$\delta(C_d,p)=\frac{(d-1)(d-2)}{2}=(r-1)(2r-3).$$
It follows that the number of branches of the singularity $(C_d,p)$ is 
$$2\delta(C_d,p) -\mu(C_d,p)+1=2,$$
and hence $C_d$ is a nearly cuspidal rational curve. Apply now Theorem \ref{thm3}
with $d'=r-1$, and conclude that $C$ is a 4-syzygy curve with exponents
$d_1=d_2=d_3=d_4=r$. 
\end{ex}

Consider now the prime decomposition of $d$, namely
\begin{equation}
\label{pfactor}
d=p_1^{k_1} \cdot p_2^{k_2} \cdots p_m^{k_m},
\end{equation}
where $p_i \ne p_j$ for any $i \ne j$.
We assume also that $m\geq 2$, the case $m=1$  being settled in Theorem \ref{thm2} (3). By changing the order of the $p_j$'s if necessary, we can and do assume that  $p_1^{k_1}>p_j^{k_j}$,
for any $2 \leq j \leq m$. Set $e_1=d/p_1^{k_1}.$

\begin{thm}
\label{thm4}
Let $C:f=0$ be a rational, nearly cuspidal plane curve of degree $d=2d'+1$,  an odd number as in \eqref{pfactor}. Then, if
$$mdr(f)=d_1 \leq r_0:=\frac{d-e_1}{2},$$ then $\nu(C) \leq 2$. In particular, if $d=3p^k$, with $p>3$ a prime number, then $C$ is  a free curve, a nearly free  curve or a plus-one generated curve with $\nu(C)=2$.
		\end{thm} 
		
		\proof Exactly the same proof as for \cite[Proposition 3.3]{DStMos} implies that
		$n(f)_j \leq 2$ for $j \geq 2d-3-r_0$. Our assumption implies
		$r=mdr(f) <d/2$ and hence case (ii) in Theorem \ref{thmPO2} is excluded.
		Moreover one has $2d-r-3 \geq 2d-r_0-3$ and one can use \cite[Theorem 3.2 (1)]{DStJump} to get
		$$\nu(C) =n(f)_{2d-r-3} \leq 2.$$
		When $d=3p^k$ and $p>3$, we set $p^k=2p'+1$. With this notation, we have
		$e_1=d/p^k=3$, $d=3(2p'+1)=2(3p'+1)+1$, and hence $d'=3p'+1$.
		Due to Theorem \ref{thm3}, we can assume $r=mdr(f)<d'$. Since
		$$r_0=(d-3)/2=3p'=d'-1$$
		the final claim in Theorem \ref{thm4} is proved.
		\endproof
		
\begin{cor}
\label{corConj2}
Let $C:f=0$ be a rational, nearly cuspidal plane curve of degree $d$. Then Conjecture \ref{c2} holds in any of the following cases.

\begin{enumerate}

\item The degree $d$ is even;

\item The degree $d$ is at most $33$.
\end{enumerate}

\end{cor}		
The last claim follows since all the odd numbers $d \leq 33$ are either powers of primes, or of the form $3p^k$, with $p>3$ prime.

\end{document}